\documentclass{article}
\usepackage[utf8]{inputenc}

\usepackage{amsmath,amsthm}
\usepackage{amsmath}
\usepackage{amssymb}
\usepackage{amsfonts}
\usepackage{amsthm}
\usepackage{amsbsy}
\usepackage{amsgen}
\usepackage{amscd}
\usepackage{amsopn}
\usepackage{amstext}
\usepackage{amsxtra}
\usepackage{mathrsfs}
\usepackage{enumitem}
\usepackage{graphicx}
\usepackage{verbatim}
\usepackage{epstopdf}
\usepackage{float}
\usepackage[all,cmtip]{xy}
\usepackage{accents}
\usepackage{sseq}
\usepackage{url}
\usepackage{makeidx}
\usepackage{wrapfig}
\usepackage{tikz-cd}
\usepackage{framed,enumitem} 
\usepackage{syllogism}
\usepackage{tikz}
\usetikzlibrary{decorations.markings,arrows}

\usetikzlibrary{matrix}
\usepackage[T1]{fontenc}
\usepackage{titlesec, blindtext, color}
\definecolor{gray75}{gray}{0.75}
\newcommand{\hsp}{\hspace{20pt}}
\titleformat{\chapter}[hang]{\Huge\bfseries}{\thechapter\hsp\textcolor{gray75}{|}\hsp}{0pt}{\Huge\bfseries}

\usepackage[margin=0.8in]{geometry}

\theoremstyle{definition}

\theoremstyle{remark}

\newtheorem*{ex*}{Ejercicio}

\newcommand{\R}{\mathbb{R}}

\usepackage{scalerel,stackengine}
\stackMath
\newcommand\reallywidehat[1]{%
\savestack{\tmpbox}{\stretchto{%
  \scaleto{%
    \scalerel*[\widthof{\ensuremath{#1}}]{\kern.1pt\mathchar"0362\kern.1pt}%
    {\rule{0ex}{\textheight}}%WIDTH-LIMITED CIRCUMFLEX
  }{\textheight}% 
}{2.4ex}}%
\stackon[-6.9pt]{#1}{\tmpbox}%
}

\usetikzlibrary{matrix,arrows,decorations.pathmorphing}
\tikzset{commutative diagrams/.cd}

\usepackage{hyperref}

\usepackage{natbib}
\usepackage{setspace}
\usepackage{afterpage}
\usepackage{sectsty}
\sectionfont{\large}

\begin{document}

\title{Wittgenstein on decisions and the mathematical practice\footnote{The first draft of the article was finished on December 5, 2022. This version is from February 6, 2023.}}
\author{M. Muñoz Pérez }
\date{\texttt{mmunozpe159@alumnes.ub.edu}}
\maketitle

%1.25
\setstretch{1.25}

\begin{abstract}
    Putnam and Finkelstein can be read as providing an answer to Kripke's skeptical argument by appealing to the way mathematics is commonly pursued. Nowadays, the debate surrounding pluralism has questioned the postulation of a unique way of developing mathematical activity. In this paper, we wish to reformulate Kripke's argument as a challenge for the conjunction of \textit{ifthenism} and a reasonable form of pluralism (which we have called $V$-\textit{pluralism}) and, at the same time, propose a reading of some passages of the \textit{Philosophical Investigations} as a solution. Our conclusion is that, in order for pluralism to be preserved, we need to clarify both the fact that we make definite \textit{decisions} and the philosophical value that they, as we argue, bear. To investigate the nature of this value is one of the further tasks want to cover in this article.
    \\

\end{abstract}

\textbf{1. Ifthenism and other inclinations of mathematicians.} It is well known how Wittgenstein regarded the philosophical work of some mathematicians of his time that attempted to provide philosophical explanations concerning their mathematical achievements or beliefs. In short: what was (and, one could say, is) dangerous about this sort of explanations is that, with them, `genuine scientific work masquerades as metaphysics, rather than the other way around'\footnote{Floyd, J. (2000) \textit{Wittgenstein, mathematics and philosophy}, p.235, in \textit{The New Wittgenstein}, eds. Crary, A. \& Read, R.}. This can be clearly seen in the \textit{Philosophical Investigations}\footnote{In what follows, PI. We also adopt the usual abbreviations for titles of Wittgenstein's works, such as TLP for \textit{Tractatus Logico-Philosophicus}. Moreover, we will indicate the which part of PI we are referring to with the aid of a roman numeral, followed by the number of the paragraph.}, where the inclinations of mathematicians are considered, not as something philosophical, but as something that philosophy needs to clarify, to \textit{treat}\footnote{PIi§254 and §255.}. Indeed, sometimes mathematicians intend to describe something they find deep or relevant and philosophical pictures come to their aid, but they may carry with them some misleading properties, concerning their own work and how it relates with other ramifications of knowledge. Note that this use (the positing of philosophical pictures) clashes with the characterization that early Wittgenstein makes of philosophy as an \textit{activity} rather than a set of `philosophical propositions' or a `theory'\footnote{TLP4.112.}. We will be revisiting this attitude in the next sections, since one of our main aims is to clarify the sense in which a danger is here present. 
\\

One could argue that this kind of problems only used to arise when the platonistic or set theoretic views received more attention, that is, when the mathematical community was not as diverse as is today and not as aware of crucial metamathematical results\footnote{Some, as Putnam, argue that philosophy of mathematics is now harder than before the publication of Gödel's theorems.}. But this claim overlooks the intricacy of some actual views held by the nowadays working mathematician. While it is true that it seems that, in some areas, being `foundations-agnostic' is perfectly fine, this situation is far from clear in topics such as category theory and set theory (and even these samples are insufficient to provide a faithful depiction of some of the latest revolutionary mathematical frameworks). It is not our aim to enter here into the discussion of the role of this disciplines in the grounding (if there is supposed to be) of mathematical knowledge. Our only concern may be that some attitudes towards philosophy -as irrelevant, redundant or superfluous- actually disguise (philosophically-charged) pictures and that this yields the same kind of problems that we have described above: mathematical results as involving items of `superlative' nature\footnote{PIi§192.}. 
\\

For example, the teacher must be able to present set theory as answering some issues concerning the infinite while, at the same time, dismiss their importance for the implicit development of the theory. One who has learned enough category theory might be inclined to hypothesize that certain abstract structures provide a `better' foundation of mathematics than sets. What this situations both have in common is that they present some mathematical results as having special importance, as grounds for accepting some \textit{striking} philosophical ideas\footnote{PIi§219. We know that Wittgenstein was against this idea of surprise in mathematics or, more precisely, of mythical explanations of established facts: cf. PIi§221.} but, at the same time, they are endangered by the potential problem of falling into metaphysical nonsense. 
\\

A further phenomenon of today mathematics is what we will be calling \textit{ifthenism}\footnote{Although we are aware of the use of this term in other debates (such as the one concerning modal structuralism) we will use it to name the view presented in the text because of its similarity with inferentialism. In fact, ifthenism could be considered a consequence of inferentialism when we restrict ourselves to different (formal, logico-mathematical) contexts.}. The core idea of ifthenism is that a mathematician acquainted with different formal calculi is not committed by his work to hold private beliefs concerning preference or grounding for preference, in the sense that he only works hypothetically: he does not believe in the axioms, for example, of ZFC, but he is certain that anyone who believes them will agree with the theorems of that theory. In other words, the mathematician is not to be concerned with the problem of devising a calculus in order to capture some insight but with the task of deriving (one may add, blindly) consequences from its axioms\footnote{Note how this position could be favored by the working mathematician who calls himself agnostic when some philosophical issues arise in his field of research. What we argue is that the working mathematician accepts some philosophical pictures but not others, being those dismissed as something alien to mathematical practice itself. Normally, this attitude is disguised by simply stating that \textit{mathematics must be left for the mathematicians}, leaving no room for useless and intrusive philosophical comments.}. Thus, ifthenism can be thought as naturally accompanied by some form of \textit{pluralism}, at least when the mathematician we are considering knows how to work in different frameworks, as is usual nowadays\footnote{The connection between both views is not that obvious. At a first glance, one could say that, if ifthenism is claimed as a fact of common mathematical practice, it would then (philosophically) enable pluralism, in the sense that there would be no tension at all in the fact that the mathematician understood properly different notions of, for example, logical consequence, and could at the same time work with all of them in different ways. We will study deeper this connection in the following.}: this would add that, in fact, we can provide different calculi for the same informal content and that all of these are equally valid in some sense. We can see how pluralism fulfils the philosophical task that ifthenism ignores from the beginning; what both of them look for is a convincing argument for the usual practice of mathematicians, which is not restricted to a fixed calculus (or collection of calculi) and therefore is independent of any corresponding monism at all\footnote{Here we can also find philosophers whose desire is to provide a serious philosophical analysis of mathematics and, at the same time, argue that mathematical practice must be taken \textit{simpliciter}. If mathematics and philosophy are to be separated in some way, we should give these analysis some \textit{status}: would it labeling them as `non-philosophical' be of some aid, if any? Of course, the \textit{status} of the previous slogans for `philosophical agnosticism' is far from being clear.}. In this article we wish to ask: how independent is the work of the mathematician from the motivations of the formal calculus in which he works? And what can Wittgenstein tell us about this?
\\

It seems clear, from what we have already seen, that Wittgenstein would favor that `no mathematical discovering can have any bearing on the philosophy of mathematics'\footnote{Dummett, M. \textit{Wittgenstein's philosophy of mathematics} (1959), The Philosophical Review, Vol. 68, No. 3 (Jul., 1959), pp. 324-348.}. Wittgenstein would also defend, at least in principle, the independence of mathematics from philosophy, as can be deduced from PIi§124. This should be enough to have any reading of PI defending some form of constructivism or finitism considered as misguided. At least in PI, it seems that there are no clear signs of open and active mathematical revisionism\footnote{Maybe the closest example the reader can think of occurs during the considerations about negation from PIi. Note that we are leaving aside comments included in RFM and CLFM.}. Despite all of this, and from our point of view, some tension arises between the mentioned paragraph and PIii§271, where psychology is equated, regarding its `[conceptual] confusion and barrenness', with set theory\footnote{Cf. TLP6.031 where it is stated that it is superfluous for mathematics.}. It is far from obvious that someone who believed in the complete separation of mathematics and philosophy would defend this kind of statements surrounding set theory\footnote{Note, however, that Wittgenstein would talk only about this link in some \textit{active} way. Philosophy would have to analyze the use of language made in set theory, at least when the mathematician makes, for example, some obscure observations regarding the general idea of some argument. (Dummett could put here the example of the diagonal argument and Wittgenstein's attitude towards it.) The working mathematician that has no time for useless debates is certainly misguided by the way in which he regards philosophy: he confuses the problem with the solution. As we will see, skepticism towards `philosophy' is usually paralleled by skepticism towards philosophical positions, as sets of beliefs, even when philosophy cannot be identified with (the adoption of) none of those.}. 
\\

\textbf{2. `Kripkenstein' and the gulf.} The famous `skeptic paradox' formulated by Kripke\footnote{Kripke, S. A. (1982) \textit{Wittgenstein on Rules and Private Language}, Cambridge, Harvard University Press.} presents a challenge for both the philosophies of mind and mathematics, as the author himself stresses in several occasions during his book. More specifically, Wright\footnote{Wright, C. (2007) \textit{Rule-following without reasons}, in \textit{Wittgenstein and Reason}, ed. John Preston; RATIO Volume XX no. 4.} sees the problem raised by Kripke as an attack to (at least) one condition that rules have to meet in order to lead us properly: `how and when can it have been settled that it is one specific rule in particular which we are following when everything we may so far have said, or explicitly thought, or done would be consistent with its being any of an indefinite number of potentially extensionally divergent rules?'\footnote{p. 482.} Kripke illustrates this with the example of the addition symbol and the intended operation `quus'. One of the skeptic arguments he gives reads as follows:
\begin{quote}
    It might be urged that the quus function is ruled out as an interpretation of `$+$' because it fails to satisfy some of the laws I accept for `$+$' (for example, it is not associative; we could have defined it so as not even to be commutative). One might even observe that, on the natural numbers, addition is the only function that satisfies certain laws that I accept- the `recursion equations' for $+$: $(x)(x+0=x)$ and $(x) (y) (x+y'=(x+y)')$ where the stroke or dash indicates successor; these equations are sometimes called a `definition' of addition. The problem is that the other signs used in these laws (the universal quantifiers, the equality sign) have been applied in only a finite number of instances, and they can be given non-standard interpretations that will fit non-standard interpretations of `$+$'. Thus for example `$(x)$' might mean for every $x<h$, where $h$ is some upper bound to the instances where universal instantiation has hitherto been applied, and similarly for equality\footnote{This appears in footnote 12, pp. 16-17.}.
\end{quote}

In his paper against Wright's answer to Kripke, Finkelstein\footnote{Finkelstein, D. H. (2000) \textit{Wittgenstein on rules and platonism}, in \textit{The New Wittgenstein}, eds. Crary, A. \& Read, R.} regards skepticism and platonism (together with the views of Wright) as misguided and alien to Wittgenstein's philosophy. The problem posed by the skeptical paradox is apparent: it only arises because we assume that there is in fact some `gulf'\footnote{As the interlocutor says in PIi§431.} between the rule and what should be done in order to satisfy it. Once we start thinking about rules as embedded in our lives, we see that there is no gulf at all\footnote{At least, that happens in all normal cases. Finkelstein refers to PIi§85.}. As Finkelstein puts it, the way we overcome the apparent problem is by means of the `weave of life'\footnote{He appeals to PIiii§2.}, not because it takes the place of interpretation or institutions in building a bridge between rule and action but because it grants that there is (normally) no need for such bridge. So the conclusion at which Finkelstein finally arrives is that, once we keep in mind how rules work in usual and well known contexts, it is allowed to talk again about the autonomy of the rule as proposed by the platonist, without adding anything mysterious to it. 
\\

The main problem that Finkelstein leaves open is that of the scope that the `weave of life' is intended to have. Are we to agree that every context we can think of falls into to the tapestry of our lives in the same way? Does this include the scientific discourse? In Finkelstein's paper one can find very few remarks explaining how his view should be extended to mathematics. He is reluctant to define the view supposedly held by Wittgenstein as `innocent platonism' (that is, platonism without any `mythical' appeal)\footnote{Check the footnote 37 of his paper.}. Nevertheless, it is possible to claim that his ideas are closely related to the reading of Wittgenstein's philosophy of mathematics made by Putnam\footnote{Putnam, H. (1997) \textit{On Wittgenstein's philosophy of mathematics}, Proceedings of the Aristotelian Society, New Series, Vol. 97, pp. 195-222.}. 
\\

\textbf{3. Putnam and commonsense realism.} After closely analyzing Kripke's arguments, Putnam concludes, similarly to Finkelstein, that the skeptic paradox was not a big deal for Wittgenstein at all, that Kripke ends up confusing Wittgenstein views with those of the interlocutor from PI\footnote{p. 259.}. Putnam certainly wants to attribute a rule-following `commonsense realism' to Wittgenstein, that is,  `not [...] the expression of scepticism of any kind about \textit{rule following}, but rather the expression of a `scepticism' about philosophical discussions of rule following'\footnote{p. 261.}. It is easy to compare this assertion to the `no-gulf' slogan of Finkelstein. But now, Putnam is ready to make the following reformulation regarding mathematical discourse:
\begin{quote}
    In short, just as there is perfectly ordinary way of learning the concepts `do so-and-so \textit{ad infinitum}', `keep doing the same thing', `do this uniformly', so there is an ordinary way of learning the use of such expressions as `either $S$ is true or $S$ is false' in mathematics; and just as I take the point of the discussion of rule following not to be the expression of scepticism of any kind about \textit{rule following}, but rather the expression of a `scepticism' about philosophical discussions of rule following, so I would urge that a Wittgensteinian attitude towards the use of the law of the excluded middle in mathematics should involve not a scepticism about the applications of that law within mathematics, but a `scepticism' about the very \textit{sense} of the `positions' in the philosophy of mathematics\footnote{p. 261.}.
\end{quote}
Here we should separate the following theses:
\begin{itemize}
    \item[(RL)] There is an ordinary way of learning the use of rule-following expressions.
    \item[(ML)] There is an ordinary way of learning the use of LEM and other common expressions in mathematics\footnote{Note that, when Putnam is talking about LEM, i.e. `either $S$ is true or $S$ is false', he is not leaving open the interpretation of `true' and `false': he \textit{is} pointing at the usual notion of truth and not, for example, that one of derivability, as the intuitionist would do.}.
    \item[(RS)] Philosophical explanations regarding rule-following are pointless. 
    \item[(MS)] Talking of philosophical positions in philosophy of mathematics is pointless.
\end{itemize}
Then, the argument takes the following form: we assume that (RL) implies (RS) and, since (ML) is (in some way) related to (RL), we infer (or it seems at least natural to say) that (ML) implies (MS)\footnote{Putnam seems to suggest that there is no reason against the idea that, if the first implication holds, the same argument can be applied to the second case.}. Early on the text, Putnam accepts that Wittgenstein could be committed to commonsense realism about truth only through commonsense realism about provability (an instance of (RL)), but that this is impossible given Wittgenstein's attitude towards the phenomenon of undecidable statements in mathematics\footnote{pp. 249-250.}. Nevertheless, Putnam quotes PIi§516 as providing a ground for arguing that a Wittgensteinian should accept his reading of the problem\footnote{He also provides some arguments based on how physical science works and the fact that any account of mathematical statements (regarded in his paper as `mixed') should include some satisfactory way of dealing with the physical. He claims that only provability, together with verificationism, is insufficient to provide such an account.}$^{,}$\footnote{It is quite interesting to note that the same paragraph cited by Putnam in order to support his view is cited by Dummett in his paper, together with some (revisionist) observations that Wittgenstein made against set theory. Dummett concludes from here that the  `Wittgensteinian' claim about not interfering with mathematics should be regarded as inconsistent and therefore as no well-formed claim at all. But the problem here is to determine the way in which this link is to be explained since, as we have stressed before, Wittgenstein seems to be far from being a radical revisionist of mathematics (of course, Dummett attributes him a finitistic philosophy of mathematics).}. 
\\

But let us consider the way in which Putnam approaches the term \textit{ordinary}. He subscribes the view that the ordinary should be distinguished in Wittgenstein from the philosophical: it includes all scientific discourse and, at the same time, in the case of mathematics, it does not only consist merely in proofs but acknowledges other usual activities such as conjecture-making or the application of theorems\footnote{p. 260.}. From here, it is easy for Putnam to note some immediate differences between both contexts: the phenomenon of undecidable statements is only problematic when approached philosophically, it is completely trivial while working in an ordinary context\footnote{But it is not really clear what Putnam tries to illustrate when referring to a `philosophical context'. He says in p. 260:
\begin{quote}
If it occurs in the context of a philosophical discussion as to whether undecided propositions are `really' true or false, then a Wittgensteinian may well doubt that it makes any sense, because of the peculiar metaphysical emphasis on the notion of being \textit{really} true or false.
\end{quote} Now, it is clear that Wittgenstein himself would argue against this kind of philosophical expressions. But Putnam's view ultimately rests on the supposition that philosophy \textit{is} like this. Wittgenstein would claim that, in fact, philosophy \textit{should} help \textit{within} the ordinary contexts too, at least when necessary. Compare this with the remark from PI on set theory we have quoted earlier.}. So Putnam favors, at least in the way expressed, some kind of division between philosophy and mathematics. We can try and ask: if mathematics is to be separated from philosophy, how could the `ordinary' be revised, that is, how can we have a (more or less) precise delimitation of what these ordinary contexts consist of? Of course, Putnam would see this question as belonging to the same class of those that he pretends to `take apart patiently'\footnote{As he says in p. 262 regarding some possible objections.} or, as it could also be, he would probably confine it to the limits of the ordinary itself, where science progresses independently. 
\\

Putnam says that the way in which Wittgenstein -in his reading- pretends to deal with problems of the philosophy of mathematics is to actually depart from the imposition of patterns of truth in a similar fashion as when working in the context of physical science, that what we should do instead is `to understand the life we lead with our concepts'\footnote{p. 263.}. It is true that this meets quite softly with the anti-systematic spirit present in the later Wittgenstein, being the analogy with the `forms of life'\footnote{PIii§345. Furthermore, one could compare this general attitude that Putnam's Wittgensteinian embraces with the corresponding passages of TLP concerning \textit{showing}.}, for example, quite clear. But it is notorious that this pernicious influence that Putnam claims physical science to have on our approach to some questions of philosophy of mathematics could be applicable for some of the same arguments that lead him to defend his conclusions\footnote{Compare this with what we said above. It is true that any proper philosophy of mathematics should acknowledge the relevance of application and applicability to science (even when this only consists in saying that this relevance is only apparent). But the kind of arguments that Putnam makes rest on the fact that we should `take our physics seriously', and that seems to mean that physics should play a relevant role in the philosophy of mathematics. That is, Putnam seems to evaluate some position of the philosophy of mathematics (constructivism or provability-based explanations) by establishing its insufficiency, together with the rejection of verificationism, as an explanation of mathematical applicability. It is then clear that he needs some extra suppositions in order to take down the philosophical positions he wishes to attack in the first place (constructivism), and that these suppositions come from some considerations surrounding physics.}. 
\\

\textbf{4. Ordinary language and pluralism.} Let us remember how Kripke challenge has been met. Putnam tells us that, essentially, we can read Wittgenstein as defending (RL) and therefore its implication (RS). This reading is, for the most part, consistent with PI and, as such, we are willing to accept it. On the other hand, nevertheless, we have seen how there are very few parts in PI (even ignoring the contradictions that arise with those in other works of Wittgenstein) where the commonsense realism vindicated by Putnam for rules applies for the case of mathematics. The transition that Putnam wants to make from one scenario to the other is quite subtle, even when considered outside the interpretative debate of what Wittgenstein (or a Wittgensteinian) should be committed to. While we take the Putnam-Finkelstein answer as legitimate, that is, as an important observation regarding how Kripke's problem cannot attack our faculty of following rules in well known cases or, equivalently, as how there is no gulf between rule and action in ordinary contexts, we claim that the relevance of Kripke's argument should not be overlooked when arranged properly. After all, the main source of disagreement with the mentioned answer is the demarcation between \textit{specialized} and \textit{ordinary} contexts of linguistic use and that amounts, again, to the problem of separation between mathematics and philosophy.
\\

Now, what we think that has been made evident by our previous observations about Putnam is that he draws a very concrete picture about mathematics and mathematical practice from his remarks on Wittgenstein's views. The following lines characterize \textit{mathematics as it is actually done}:
\begin{quote}
    In the sense of having a translation into ordinary language, of course mathematics \textit{has} a semantics -a `realist semantics' if you please. Or rather, mathematics as it is actually done, `unregimented' mathematics, already \textit{is} in ordinary language\footnote{p. 262.}.
\end{quote}
Of course, the adjective `unregimented' can be read as pointing against a constructive account of mathematics\footnote{For example, pp. 249-250 contain some direct references to this view. This, of course, is unsurprising, since Putnam's paper is based on the idea of reading Wittgenstein at not necessarily committed to the whole of the `antirealist story'.}. It is not difficult to find this point controversial\footnote{We could think of further issues concerning Putnam's `realist semantics'. For example, one could argue that Wittgenstein is attracted to embrace in some way the context principle (supposing, of course, that this is in tension with realist semantics): see for example PIi§49 and, above all, PIi§544:
\begin{quote}
    When longing makes me exclaim “Oh, if only he’d come!”, the feeling gives the words ‘meaning’. But does it give the individual words their meanings? 
    
    But here one could also say that the feeling gave the words \textit{truth}. And now you see how the concepts here shade into one another. (This recalls the question: what is the \textit{sense} of a mathematical proposition?)
\end{quote} All this, of course, leaving aside the famous remarks from TLP3.3, 3.314.}. Do all mathematicians work in a well defined, definite and absolute way? Is this picture sound? What Putnam says can be easily contrasted with the main ideas that mathematical pluralism wishes to defend. If we follow Zalta\footnote{Zalta, E. N. \textit{Mathematical Pluralism}, the Metaphysics Research Lab.}, it seems reasonable to argue that realism has been challenged -as the natural philosophy for the working mathematician- by pluralism. In short: there are \textit{several} ways in which the working mathematician develops her activities (at least, all those included by Putnam in the ordinary language of mathematics) and there is no need of `reducing' one theory to another in order to apply some kind of indispensability argument. Therefore, the debate concerning commonsense realism of Putnam leads us to take a closer look of the pluralist mindset. 
\\

\textbf{5. Logical pluralism and $V$-pluralism.}  Beall and Restall\footnote{Beall, J.C. \& Restall, G. (2000) \textit{Logical pluralism}, Australasian Journal of Philosophy 78 (4), pp. 475 – 493.} present \textit{logical pluralism} as the adoption of the following assertions:
\begin{itemize}
    \item[(L1)] The \textit{intuitive notion} of (logical) consequence is given by the principle
    \begin{quote}
        $(V)$ $\varphi$ follows from $\Sigma$ iff in any \textit{case} in which each premise in $\Sigma$ is true, is also a case in which $\varphi$ is true.
    \end{quote}
    \item[(L2)] A \textit{logic} is given by some `\textit{specification} of the cases' that appear in $(V)$\footnote{In p. 2 we can read: `Such a specification of cases can be seen as a way of spelling out \textit{truth conditions} of the claims expressible in the language in question'. That is, they would include the usual conventions regarding the primitive symbols of the logic.}.
    \item[(L3)] There are at least two different specifications of the cases that appear in $(V)$.
\end{itemize}

Now, once one takes this picture for granted, it seems natural to define a mathematical framework as simply what mathematics consist of when `pursued in the context of' some logic, as defined in (L2)\footnote{pp. 11-12. We can also read: `Theorems of constructive mathematics are simply theorems of mathematics proved constructively. [...] the theorems of constructive mathematics are also theorems of classical mathematics. The difference between constructive and classical mathematics is not one of subject matter, but one of the required standards of proof.'}. Here, one would say, the contrast between unregimented and constructive mathematics is made in virtue of a difference of proceeding differently in each case, according to a different logic. Hence, we would be allowed, in principle, to consider mathematical statements \textit{simpliciter}, while the methods of proof would be restricted to some instance of $(V)$. 
\\

The case of logical pluralism allows us to see, nevertheless, a picture which is common among the (wider) position of mathematical pluralism. There is no reason why, we argue, we should not consider other intuitive notions that admit several formal specifications. Namely, consider the following scheme as an extension of the one given above:
\begin{itemize}
    \item[(V1)] The \textit{intuitive notion} of a (mathematical) concept $\chi$ is given by some principle $(V_\chi)$ in which some \textit{cases} occur\footnote{Is the notion of logical consequence a mathematical one? Well, if we define mathematics as the study of formal systems, is clear that by definition the concept $\chi$ would be mathematical, in the sense that it is studied through the formal properties of a system generated by the specifications of some case from $V_\chi$.}.
    \item[(V2)] A \textit{formal system} is given by some \textit{specification} of the cases that appear in $(V_\chi)$. 
    \item[(V3)] There are at least two different specifications of the cases that appear in $(V_\chi)$\footnote{It seems, from this theses, that natural language actually plays a definite role within the formalization process of $\chi$. One would like to read here some of the insight from TLP6.233.}.
\end{itemize}

Let us call the adoption of (V1)-(V3) (or rather, of their instances) $V$\textit{-pluralism}. As we foresaw, $V$-pluralism is based in a picture which is quite distant from that defended by Putnam's Wittgensteinian, for one could consider what is included in the domain of the ordinary language of mathematics and still view it as susceptible of further specification (in the sense above), being unclear at all its uniqueness or, alternatively, one could simply argue that there are more possible specifications to the contents supposedly formalized by this ordinary language. Note that $V$-pluralism might defend that it is part of the mathematician's work to actually provide these specifications and therefore Putnam would be lead to admit that this usual activity would be alien to ordinary mathematical practice. To put it briefly, while $V$-pluralism would admit with Putnam's Wittgensteinian that there is a common way of referring to, say, mathematical conjectures, this way of talking about \textit{truth} would be in fact embedded in a narrower notion dependent on the formal system one is working with (not, of course, merely consisting in the `unregimented' case).
\\

It is quite remarkable that $V$-pluralism gives us a position in which two different kinds of content, `intuitive' and `formal', are linked. That is, in a situation as described by (V1) we would implicitly have all the problems of a \textit{philosophical discussion}: how is the principle $V$, for instance, decided as \textit{the} expression of the initial insight? In situation (V2) we would similarly have a process, namely the \textit{formalization process}, that begins with some \textit{reading} of $V$ and ends with the (basic specifications of the) formal calculus. One form of ifthenism would consist in simply saying that the mathematician needs only to be concerned with the study of the properties of the formal system that arises from (V2). Another view might say that the mathematician needs to be equally concerned with the basic insight from (V1). What is clear is that $V$-pluralism implicitly holds the following:
\begin{itemize}
    \item[(1)] The mathematician works following some initial intuition or insight that formalizes by means of some calculus. 
\end{itemize}

Now, the precise relationship between logical and $V$-pluralism depends essentially on the role we assign to logic with respect to mathematical statements. One could say (as we pointed out above) that mathematical statements could be considered \textit{simpliciter}, but that what actually mattered would be the fact that we could prove them within some logic, or one could be persuaded to reply that the mathematical content would be modified together with the logic used by the corresponding formal system\footnote{Compare this with what we said above. Note that one could read Putnam (\textit{contra} what we have said earlier) as not entirely going against constructivism and just stating that there is a way in which we ordinarily understand mathematical statements such as definitions or conjectures. In this case, the $V$-pluralist might argue that Putnam is simply accepting the existence of statement $V_\chi$, but note that -from his point of view- a further distinction is needed by means of some cases specification from such $V_\chi$. Moreover, it is not obvious how the $V$-pluralist would regard the conception that Putnam draws consisting in philosophical debate as opposed, v.g., to mathematical discourse.}. Take the following example. The question would be: what is the nature of the statement `$f: \R \to \R$ is continuous or is not'? Several answers could be given by a pluralist:
\begin{itemize}
    \item[(a)] We can always understand the statement `$f: \R \to \R$ is continuous or is not', but there are many ways in which we could formalize `$S$ or not $S$'.
    \item[(b)] We can always understand the statement `$f: \R \to \R$ is continuous or is not', but there are many ways in which we could formalize `$f: \R \to \R$ is continuous'\footnote{It would seem at least reasonable to say that Putnam's Wittgensteinian could defend a similar claim to this: remember that he was particularly skeptical about understanding `$S$ or not $S$' in a different way than the standard one. Against this, note that Putnam seems to say that, because we can understand the statement, there must be a concrete specification that formalizes it. Nevertheless, Putnam can be read as arguing, in general, that classical mathematics are already given in ordinary language, so it also seems that he would stand against the idea of providing different formalizations of `$f: \R \to \R$ is continuous' (he does not concrete examples, but the paragraph from PI he quotes can be read as implying this). Therefore, in this case we could only attribute him the view that we commonly understand such statements and thus, this would be an incomplete position for the $V$-pluralist (see above).}.
    \item[(c)] We can always understand the statement `$f: \R \to \R$ is continuous or is not', but there are many ways in which we could formalize the statement altogether\footnote{Wouldn't the working mathematician, then, fall into philosophical misfortunes when trying to talk in a ordinary context? Not at all, for it would be clear at every moment in which formal system he is thinking about. At least, since the advances of metamathematics, it should be more natural to think that one formal system is different from another in the same way we say that a given function is different to another.}.
    \item[(d)] We only understand `$f: \R \to \R$ is continuous or is not' when we are talking within a concrete context, given by some specification concerning `$f$ is continuous' or `$S$ or not $S$'.
\end{itemize}
We could name the positions behind (a) and (b), respectively, \textit{pure logical pluralism} and \textit{pure mathematical pluralism}. Now, it seems natural to argue that $V$-pluralism would indeed give an answer similar to (c). Of course, (d) could not be held by the $V$-pluralist, since thesis (V1) requires an informal formulation that \textit{later} leads to some specification, not reciprocally\footnote{That is, a view as the one expressed by (d) would only make sense if read as saying: `we can only understand \textit{formally} the statement from the standpoint of a given formal calculus'. This is, in the context we are discussing now, vacuously true.}. Therefore, $V$-pluralism could state that there are different kinds of mathematical concepts formalized within a fixed calculus, each fitting into a precise role\footnote{Or, more precisely, one could say that a formal calculus is designed to study a set of notions. An example of this situation would be that of ZFC assuming the (classical) first-order logic framework. Here, $V$-pluralism would simply see the assumptions of first-order logic being made explicit together with those of ZFC.}. 
\\

Now, let us link $V$-pluralism to two common formulations of mathematical pluralism. According to Zalta\footnote{See Zalta's paper above.}, there are three main forms in which mathematical pluralism can be presented:
\begin{itemize}
    \item[(A)] Every mathematical theory which is consistent consists of truths about a domain of objects associated to it\footnote{Although we will not be entering in ontological issues, the notion of a domain of objects \textit{associated to} a theory can be linked with Quine's ontological commitment. Nevertheless, note that, as Zalta says, `[a] pluralist leaves the determination of what is interesting to mathematical
    practice'. Also, mathematical pluralism finds natural to question the idea of reducing one theory to another, as if, in that case, their domains of objects would have to be the same.}.
    \item[(B)] Every mathematical theory, consistent or inconsistent, consists of truths about a domain of objects associated to it.
    \item[(C)] Each position in philosophy of mathematics is based upon an insight about the nature of mathematics that is valid in some sense\footnote{Note that the contrast with Putnam (or his Wittgensteinian) becomes even more evident in virtue this thesis.}.
\end{itemize}
It seems natural to argue that the `insight' that appears in (C) could be related (if not the same as) the one appearing in (V1). For example, some revisionist account of mathematical nature would probably lead to the revision of the basic intuitions that some (classical) mathematical theory could have already formalized\footnote{Here, it is true, we could equally claim that the intuition of (V1) would be the same for everyone but that different specializations could be given or that, differently, each of these specializations rested on different intuitions. That is, logical and $V$-pluralism, as presented, have to deal with the fact that we should be able to grasp different things from the same statement $V$. If we wished to restrict our abilities to detect different `cases', much of the appeal of pluralism would be lost. Nevertheless, Putnam could always argue that, in order to formulate pluralism itself, we must be `non-ambiguous' and that therefore we assume some definite logic we intend to other people to follow (in other words, since we are in meta-discourse, we should agree on some meta-logic as background for understanding). This actually corresponds to the `one-many' answer that can be read in p.17 of Beall, J. C. \& Restall, G. \textit{Defending Logical Pluralism}, in \textit{Logical Consequence: Rival Approaches}, eds. Brown B. \& Woods, J., pp. 1-22: `Graham Priest poses the question: “Logic: One or Many?” Our answer is “both”.
\begin{quote}
\textit{One}: There is precisely one core notion of logical consequence, and
that notion is captured in schema (V).

\textit{Many}: There are many true instances of (V), each of which specifies
a different consequence relation governing our language.   
\end{quote}
This one-many answer is what we call ‘pluralism’.' We can overlook the difficulties of such objection for the moment, but note that we could simply make the $V$-pluralist to adopt some form of scepticism regarding this background, underlying logic.}. So, as a consequence of this view, $V$-pluralism could be viewed as belonging to kind (C) of mathematical pluralism. From this, we could derive the following reformulation of (1):
\begin{itemize}
    \item[(1')] The mathematician works following some initial (philosophical) insight about the nature of mathematics that formalizes by means of some calculus. 
\end{itemize}

In addition, provided that we were able to encompass, generally, a definite domain of objects with the aid of the different cases that appear in the principle $V$, this would mean that the specialization of one of these cases would lead to the theory that has that domain as the associated one. That is, we could also see $V$-pluralism as stating something similar to (A). Now, it may be argued that a proper mathematical pluralism is to be committed to assert that any given formal system\footnote{This is what Zalta says about kind (A) of pluralism. Note that we are talking about \textit{consistent} formal systems so that we restrict ourselves to this case.} is equally valid, no matter the arbitrariness of its axioms, since this is a fact only revealed by later mathematical practice\footnote{In p. 2 of  Zalta's paper it can be read: `A pluralist leaves the determination of what is interesting to mathematical practice', although `not every theory is equally fruitful'.}. This would imply, in principle, that $V$-pluralism should subscribe this view: it could be that some formal system could bear some application in the future, the intuition of that system would then derive from this application and, thus, this formal system would satisfy the requirements imposed by (V1)-(V3)\footnote{This is a pluralist account of mathematical \textit{applicability}: if there is a useful application of a formal system, it then seems natural to take this fact itself as something `informal', in the sense of `susceptible of formulation within natural language in some statement $V$'. In this case, the schema (V1)-(V3) could be repeated again.}. Thus, a $V$-pluralist will defend some form of the following assertion:
\begin{itemize}
    \item[(2)] Every formal modification of the calculus yields an equally valid or acceptable one. 
\end{itemize}

A natural question, regarding what we said at the beginning, may arise: what is the relationship between $V$-pluralism and ifthenism?
\\

\textbf{6. Kripke's \textit{skeptic paradox} revisited.} Take a look at Kripke's text quoted above. Remember that Wright's formulation of the skeptic paradox was that, up to some point, we had only a finite set of instances were we applied a given rule, so instead of the rule we were pretending to follow there was no reason to belief that we were actually following another (which differed extensively with the other, the first deviation being in a case at which we had not yet arrived). The Putnam-Finkelstein answer showed how this challenge was no challenge at all. As we announced, we consider this solution as legitimate when the skeptic is questioning our well-established abilities of rule-following. The problem may arise in, as we put it then, specialized contexts\footnote{Let us explain this briefly. Finkelstein, \textit{contra} Wright, tries to show that nothing special happens when we follow a recipe, for example, and that this excludes the need for making some previous `decision', in Wright's sense. Now, we could think that following a recipe is related with its finality in a different way from how doing mathematics, as mathematical practice, is. This is what we have tried to point out with the use of `specialized' and `ordinary' contexts. One could go on and argue that, while a deviation of a recipe may affect the desired result, a modification in, say, a set of rules of a game, leads to a new game (compare this with the equally valid formal systems of (2)). This idea can be found in Valdés' introduction to the Spanish edition of Wittgenstein's \textit{Phenomenology}. Now, we should make clear whether the pluralist could be committed to such a view or not, but this at least shows that Finkelstein's view cannot be (directly) extended to every context in general.}. Our claim is that a similar paradox could be reproduced in order to attack a specific form of pluralism, namely, the conjunction of both $V$-pluralism and ifthenism. 
\\

As we said above, $V$-pluralism recognizes the formalization process of some informal notion as an important phenomenon in mathematical practice. One may be inclined to say that, once all the formalization (i.e. the determination of the related specifications) comes to an end, the mathematician should leave aside any other kind of considerations and his work should consist in studying the formal derivations possible in the resulting formal calculus. This separation is what ifthenism, as we understand it here, wishes to defend. The question we make  ourselves is: how would another option look like? 
\\

On the other hand, it is clear that as mathematicians are allowed to make some conjectures and propose definitions while working in a formal system. These elements certainly modify the formal calculus, or its links with the original insight given by the thesis $V$. Now, suppose a mathematician gives some formalization of a concept $\chi$ by specifying a case of $V_\chi$. It is clear that, up to that moment, she has only made finite decisions on the definitions and the general purposes and aims of the system. The question is, then: how is she to determine if the formal calculus she is working in is the specialization of the concept she intends to formalize and study thereafter? Is there a guiding intuition behind all the decisions she actually makes when studying the mentioned concept? With words that Kripke used in the quoted paragraph: how does she know that she is not developing a `non-standard' calculus?\footnote{An immediate answer could be: how are you talking about `intuition' or `insight'? Are not this psychologistic concepts that one should undermine instead of making some use of them? Well, if we take a look at (P1)-(P3), it is clear that this answer should attack also logical pluralism itself and this is a further debate we are not concerned about.} This situation gets even more clear if we admit for some specifications to be made \textit{while} the formal system is being developed: think about how one could add axioms as a result of some undecibility result, for instance.
\\

But let us take a closer look. We have seen how $V$-pluralism implies (1) and (2) above. Note that ifthenism will defend the following thesis:
\begin{itemize}
    \item[(3)] The mathematician needs no intuition or insight for developing his work. 
\end{itemize}
This clearly goes against the crude form of $V$-pluralism as expressed by (1). So the conjunction of both thesis would in fact give us:
\begin{itemize}
    \item[(4)] The only intuition required for mathematical practice is the one formalized by the initial specifications from which the formal system arises.
\end{itemize}
This is the affirmation we claim as depicting an incorrect picture of mathematical practice, for consider the case in which it becomes natural for the mathematician to accept some axiom (we will assume that both this axiom and its negation are undecidable from the original system alone). Both theories are, by (2), equally acceptable in some sense. What could the mathematician argue in favor of accepting one between the axiom or its negation? We could say, following (4), that it is the original intuition what makes him accept one or the other, as something that imposes itself upon anyone that understands properly the formal system. But, how does he know which, between -again- two equally valid formalizations, the one that corresponds to the initial intuition (or the one that he \textit{intends} to give)? We could somehow argue that the original intuition led to this and that situations, and that we made the correct choices: this is what (2) actually pretends to leave as open, since the value of some theory could be granted with future applications. So, finally, what we arrive at is a stage in which ifthenism, as expressed by view (4), is in tension with (1).
\\

Note the following subtlety: we are not questioning whether the mathematician actually knows what she means with an informal conjecture\footnote{See, for example, the answer (c) that, as we argued, the $V$-pluralist should give to the question of the relationship between logical and $V$-pluralism.} but the fact (assumed by ifthenism) that she does not make any relevant decision about the formal system that exceeds proving theorems or blindly deriving consequences from already given axioms and rules\footnote{In other words, if the Putnam-Finkelstein answer is to be understood as favoring some kind of \textit{rule autonomy}, what we are trying to attack is the \textit{autonomy of intuition} as a guiding item in mathematical practice or the usual attitude of separating mathematics and philosophy while at the same time allowing the making of explicit choices which are philosophically charged.}. What we are trying to introduce in mathematical practice is our decision-making while working with a given calculus, which relates, again, to the supposedly natural separation between ordinary and specialized contexts or, alternatively, between mathematics and philosophy. Again: how could ifthenism explain the arbitrary expansions of a system? Well, $V$-pluralism could try to help by admitting the successive expansions as \textit{valid} in some sense. The question arises, once again, about the passive role that intuition supposedly plays here: indeed, the resulting system is valid but how do you know that you are faithfully following your original intentions?
\\

But all that the skeptic paradox gives us in this restricted context is a different picture of how $V$-pluralism may have to deal when linking informal and formal contents. The point where the alleged paradox vanishes is this one: \textit{we actually make definite decisions}. So what remains is to change our attitude towards how we make these choices and their relationship with usual mathematical practice\footnote{Or, with other words, we have to let philosophical activity treat this decisions as mathematical inclinations.}. The picture that ifthenism defends is that mathematicians do not have to actively care about the original notions formalized by their work: it is a picture on how intuition is an autonomous force of nature that compels us to \textit{accept} some decisions \textit{without discussion} and how it remains intact during the development of a theory. Kripke's challenge applied to this case shows the need of taking additional \textit{decisions}, one may say, in the sense of \textit{positions}\footnote{Putnam does have a point when stating that his Wittgensteinian should be skeptical concerning positions in philosophy of mathematics. If we understand `position' as a set of beliefs based around obscure reasoning, it is clear that Wittgenstein would ban its use (cf. PIi\S401). Now, what we wish to defend is that there is some contact from mathematics and philosophy in the way we take (and argue about) our decision-making. Here, a position may be just a way of making some mathematical idea clear. Of course, the debate of whether some idea is problematic is reduced by the $V$-pluralist to the fact that one can be able to specify different calculi from one (same?) idea and, moreover, that this idea corresponds to different insights. So the way we should understand this is, as we announced before, by considering philosophy as an activity that takes place within mathematical practice, that does not posit anything and only tries to make explicit the way we (naturally) act.}. Whether we call them mathematical or philosophical, it is a further debate\footnote{Note that we could consider this decision-taking as something implicitly mathematical, as an \textit{act}, but that one could only \textit{show} from a philosophical point of view. We will later also make reference to this.}; what is clear is that ifthenism, following its desire to keep mathematics and philosophy completely separated, enforces an image that leads into some serious philosophical problems. To put it differently: Finkelstein and Putnam arrive at the same point as we do, that is, there is no reason to believe that we do not know how to take decisions when working with a formal calculus. But, \textit{contra} Putnam, we actually defend that philosophical problems have some kind of link with mathematical practice (again, we wish to name them in other way if necessary). Believing that philosophy is to be completely separated from mathematical practice is to actually see philosophy as inherently positing theses that interfere with mathematics (this goes against taking philosophy as \textit{activity}) and, at the same time, reinforces the mathematician's belief that philosophy is superfluous and that philosophical matters can be left aside without any (possibly dangerous) consequence at all. 
\\

\textbf{7. Wittgenstein on decisions.} In this final part we argue that it is possible to make a reading of Wittgenstein's philosophy of mathematics in PI that can be regarded as an answer to the form of the skeptic paradox presented above. More explicitly, we argue that Wittgenstein can be read as defending the need of taking decisions within mathematical practice, since intuition alone is an insufficient guide for the development of a formal system in order to fulfil its purposes. But first, we have to take a closer look at other lectures and possible objections against this.
\\

The standard account of Wittgenstein's philosophy of mathematics is, nevertheless, sometimes linked to his philosophy of language, and therefore is labeled as some form of \textit{inferentialism}\footnote{See, for example, Zalta's paper above.}: if we accept the slogan \textit{meaning is use}\footnote{Cf. the beginning of PIi, for example, PIi§43.} and we only consider meaning as the role that one element plays in a language game\footnote{Again, this idea is widely repeated among PI: for instance, consider PIi§182.}, according to some rules, it seems reasonable to argue that, similarly, what we should only find interesting in a formal system is precisely the sets of axioms and rules we work with, that is, we could make natural an adoption of ifthenism. The main detail overlooked by this inference is that the fact that we should pay attention to the use of expressions does not entail that we should only care about the syntactical manipulations made within the system; as we said earlier, we should also recognize other activities (conjecture-making, axiom-adding, etc.) as essential elements of mathematical practice\footnote{Note that the paragraph quoted by Putnam, PIi§516, could be read as saying precisely this: 'It seems clear that we understand the meaning of the question “Does the sequence 7777 occur in the development of $\pi$?” It is an English sentence; it can be shown what it means for 415 to occur in the development of $\pi$; and similar things. Well, our understanding of that question reaches just so far, one may say, as such explanations reach'.
What we say is that, indeed, we make a definite use of natural language when making conjectures, explaining certain notions, etc. A pluralist would like to make a further step in making those notions clear and study the connections they should bear within some formal system.}.
\\

The first observations about mathematics that Wittgenstein makes in PI can be found in PIi§124-129. Let us discuss briefly this set of paragraphs. Earlier, we have mentioned PIi§124 and how it contradicted -in some way- the remarks from PIii§271. PIi§125 is quite illuminating. There, one can read that the only task that philosophy has with respect to mathematics is not `to resolve a contradiction by means of a mathematical or logico-mathematical discovery, but to render surveyable the state of mathematics that troubles us a the state of affairs \textit{before} the contradiction is resolved.' So the situation is that one of discovering a contradiction in a formal system: we formalized some intuition and now a difficulty within it arises; philosophy should help us to understand this phenomenon better before taking any further step in its solution. More precisely:
\begin{quote}
   Here the fundamental fact is that we lay down rules, a technique, for playing a game, and that then, when we follow the rules, things don’t turn out as we had assumed. So that we are, as it were, entangled in our own rules.

This entanglement in our rules is what we want to understand: that is, to survey.

[...]

The civic status of a contradiction, or its status in civic life - that is the philosophical problem\footnote{Although it does not belong to the scope of the present article, it would be interesting to study the connection of this emphasis put on the contradiction and the passages in RFM and LFM usually considered akin to some form of paraconsistency.}.
\end{quote}
In other words, we suffer an \textit{unexpected deviation} from our original intentions, we now see that there is some difficulty within the specification we have been working with. As Wittgenstein says, what philosophy should help us to study deeper is the connection with our informal insight and the formal items that we have blindly followed and resulted in some state of irrevocable tension with the former\footnote{In PIii§127-128 we can see the wittgensteinian attitude towards philosophy (clear in TLP) and its active character, against `advancing thesis'.}. This `entanglement in our rules' reminds us about we have said earlier concerning the arbitrariness of our decision-making from the joint premises of $V$-pluralism and ifthenism. Wittgenstein further links this problem with the one about rule-following and intentionality. Note that, from what we have been arguing, the way one should understand rules here is different from the ordinary contexts or, at least, the cases that Wittgenstein presents with basic arithmetic can be extended to more \textit{specialized} ones, where the `common use' is not that clear\footnote{It could even be said that what we learn under the label of `usual mathematics' deals with items of knowledge that are not susceptible of discussion because of their simplicity. It is clear that pluralism is, as we remarked earlier, a philosophy for the working mathematician, not the layman. The main difference here is that the mathematician should be prepared to arrive at results which go, in every aspect, against naive assumptions of any kind. Wittgenstein appears to concede this when he talks about the role of mathematics in what we can imagine (PIi\S334, 516, 517). Of course, the realist might say: `but does this not prove that there is, after all, something \textit{true}?' It could be so, but this is not at all our problem; here we could even say that realism tries to force us to accept some philosophical position that exceeds our conclusions: what we are mainly interested in is how, when becoming acquainted with mathematical reasoning, we are able to regard a picture of mathematical activity as conceived in (1') is plausible.}. 
\\

There are several parts that can be read as if Wittgenstein were struggling with the idea of `guiding intuition' of a formal system. For example, take PIi§213:
\begin{quote}
    “But this initial segment of a series could obviously be variously interpreted (for example, by means of algebraic expressions), so you must first have chosen one such interpretation.” - Not at all! A doubt was possible in certain circumstances. But that is not to say that I did doubt, or even could doubt. (What is to be said about the psychological ‘atmosphere’ of a process is connected with that.) Only intuition could have removed this doubt? - If intuition is an inner voice - how do I know how I am to follow it? And how do I know that it doesn’t mislead me? For if it can guide me right, it can also guide me wrong.
    
    ( (Intuition an unnecessary evasion.) )\footnote{Kripke would read this paragraph as a sign of the insufficiency of intuition to solve the skeptic paradox. Finkelstein and Putnam would see here some skepticism towards the demand imposed by that paradox, namely, a philosophical account of rule following (`[b]ut that is not to say that I did doubt, or even could doubt').}
\end{quote}
This can be read as the insufficiency of applying an intuitive notion to different cases, that is, to a situation in which we have to make a \textit{new} decision. Compare with PIi§214:
\begin{quote}
    If an intuition is necessary for continuing the series 1 2 3 4 . . . , then also for continuing the series 2 2 2 2 . . .
\end{quote}
This amounts, in our reading, to the recognizing of other systems as equally valid, as $V$-pluralism does\footnote{We will not argue anything similar to the thesis that Wittgenstein was a pluralist at heart. What could be further studied is that of the role of language games, as wholes, in giving a pluralist depiction of the use of language. From this, pluralism about formal systems -as language games- would follow. Moreover, if we compare with what we said above, we could read Wittgenstein as conceiving contexts where rule-attribution is not possible at all. Hence, it would be possible to separate the game (played according some set of rules) and the process by means of which we arrive at a set of rules (note the similarity of this picture with that depicted by the pluralist). Furthermore, this appeals to an idea that widely repeats in PI, namely, the fact that, sometimes, we \textit{act} and then we try to \textit{explain}. The fact is that we act, as Putnam says, but there is also a need for explanation, debate or comparison of possible decisions (i.e. positions): we have to make our decisions explicit.}. Again, we can link this remarks with PIi§214:
\begin{quote}
    “All the steps are really already taken” means: I no longer have any choice. The rule, once stamped with a particular meaning, traces the lines along which it is to be followed through the whole of space. —– But if something of this sort really were the case, how would it help me? No; my description made sense only if it was to be understood symbolically. - I should say: \textit{This is how it strikes me}.
    
    When I follow the rule, I do not choose.

I follow the rule \textit{blindly}.
\end{quote}
Here, the italics in the last sentence can be seen as a reference to the first lines: to say that I blindly follow the rule is nothing else than a symbolical representation of a fact. So ifthenism actually appeals to some symbolic image, nothing else; it does not deal with the problem that arises when considered together with some forms of pluralism, it does not explain or helps to analyze anything. Putnam reads this paragraph as what could be the beginning of a misleading philosophical debate around rule-following: he would dismiss ifthenism on this basis, at least in our formulation, since it enforces the view that we follow some intuition as the meaning that Wittgenstein presents in this paragraph. But here is also the point at which Putnam stops any philosophical debate and places rule-following as something we are acquainted with. The problem of this approach is, as we said, to extend it to the mathematical setting without some subtle care.
\\

Now, PIi§186 is essential for grounding our reading of PI. Our previous use (and discussion) of the term `decision' is inspired by its appearance here\footnote{It could be said that our use of the term `decision' is near the one that Wright makes and against which Finkelstein argues. But let us remember that, as Finkelstein puts it, Wright is concerned with the thesis that we make a definite decision when interpreting a rule, and that this decision actually makes rule-following (as questioned by Kripke's argument) possible. Let us remember that we do not argue such a thing about how an account of rule-following should be. What we are trying to make appealing is some kind of scepticism towards a rigid separation between mathematics and philosophy.

On the other hand, Floyd (see above) makes a use of the term which we consider to be the one we are defending here. In p. 247 one can read: `[...] there is no absolute requirement -mathematical or otherwise- that we restrict the conditions of a question in the way we do. In the trisection case, it is the \textit{decision} to require that solutions be given within a particular setting, and that solutions take a particular form and be generally applicable which generates the unsolvable -hence, provable resolvable- problem'. It seems to us that the use of the term here is not a mere grammatical coincidence.}:
\begin{quote}
    “What you are saying, then, comes to this: a new insight - intuition - is needed at every step to carry out the order ‘$+n$’ correctly.” - To carry it out correctly! How is it decided what is the right step to take at any particular point? - “The right step is the one that is in accordance with the order - as it was \textit{meant}.” - So when you gave the order “+ 2”, you meant that he was to write 1002 after 1000 - and did you then also mean that he should write 1868 after 1866, and 100036 after 100034, and so on a an infinite number of such sentences? - “No; what I meant was, that he should write the next but one number after \textit{every} number that he wrote; and from this, stage by stage, all those sentences follow.” - But that is just what is in question: what, at any stage, does follow from that sentence. Or, again, what at any stage we are to call “being in accordance” with it (and with how you then \textit{meant} it - whatever your meaning it may have consisted in). It would almost be more correct to say, not that an intuition was needed at every point, but that a new decision was needed at every point.
\end{quote}
That is, we can read Wittgenstein as, in some way, going against the picture of a guiding intuition, extending this depiction until it is impossible to hold: indeed, he seems to find misguided the answer involving an insight `at every point'. Kripke's argument above showed that, from the point of view of ifthenism, we would have to extend arbitrarily the original intuition, justifying its use in many different possible situations. The way Wittgenstein uses `intuition' here reminds us of seeing something already imposed externally (even if it is an `inner voice'), \textit{contra} the choice of the term `decision', which seems to appeal to some active character of the mathematical practice, as we said above. In fact, we can find the following in PIii§349: `Of course, in one sense, mathematics is a body of knowledge, but still it is also an \textit{activity}'\footnote{How should one understand the use of the term \textit{activity} in Wittgenstein? In PIi§23 it appears to be equivalent to a form of life (`The word “language-game” is used here to emphasize the fact that the speaking of language is part of an activity, or of a form of life') where in TLP4.112 one can find it (\textit{Tätigkeit}) as essentially different from `theory' or `set of [...] propositions'. Compare with what we said above.}. 
\\

Here, we could clearly see how the core idea of $V$-pluralism, as related with what is said in (1'), would stand: the activity of mathematics would, in some way, be related with that of philosophy: since we can start different mathematical theories based on different accounts of the mathematical nature of some informal content, it can be naturally held that there is a close connection between both activities\footnote{Namely, the one that we pointed out before: that the account of some mathematical notion made by some position leads to studying some formal system.}. Moreover, this sheds light on what a fruitful connection between philosophy and mathematics would be, against the idea of separating both and therefore eliminating some means to avoid the pernicious influence of philosophical images (say, of `philosophical theses') upon mathematical practice\footnote{That is not to say that there is a pernicious influence of mathematics upon philosophy. When doing philosophy, mathematicians can be as good or as bad as philosophers when trying to clarify mathematical practice. The mistake of the philosopher may be in arguing against usual mathematical practice based on a more or less sophisticated system of beliefs, since he assumes that mathematics are composed of assertions, in the sense of theses (one could argue here that mathematical statements have no ontological charge but, as we said, we have no space left for this). On the other hand, when a mathematician feels the need of some philosophical explanation, she produces something she \textit{does} want to take at face value as philosophical. The idea that one could interpret some theorem as establishing something deeply illuminating goes against the vision of mathematics held by Wittgenstein (compare above with the mythical images or the striking of some pictures).}. Compare this with what we said about dealing with decision-making and stressing clearly what this decisions amount to in mathematics. It is not a relationship of dependence but that of mutual help; discussions surrounding decision-making do not invalidate the formal system itself in any sense, they only go against the idea that it imposes some kind of truth upon us.
\\

\textbf{8. Concluding remarks.} Remember that one wittgensteinian insight adopted by Putnam was that evaluating philosophical ideas (concerning mathematics) in the same fashion as we evaluate physical progress carries an implicit danger. Putnam's alternative consisted in `trying to understand the life we lead with our concepts in each of' the considered areas\footnote{p. 263.}: we are inclined to fall into philosophical problems because we are misguided by some preconceptions regarding truth and objectivity. This, of course, can be linked with Putnam's slogan: \textit{mathematics do not require justifications but theorems}\footnote{Putnam, H. (2000) \textit{Rethinking mathematical necessity}, p. 229, in \textit{The New Wittgenstein}, eds. Crary, A. \& Read, R. The core idea of Putnam's paper is to argue against Quine's principle of `nothing is safe from revisability' following the tradition of Kant, Frege and Wittgenstein on the \textit{status} of logical truths.}. This statement could be understood in many ways, but it is clear that it goes against the attempt of providing a \textit{philosophical} justification of mathematics. Mathematics should not depend on some miraculous derivation within a philosophical system of beliefs. Mathematical activity is, by itself, simple in that sense. But we have argued that, although we are not to doubt about mathematics and its results or applications simply because we have no further, let us say metaphysical, explanation, we should avoid, by some similar reason, a careless use of philosophical pictures. The alternative we have tried to provide is one that can, again, call to mind multiple situations. 
\\

What Finkelstein and Putnam argue amounts to stating that we all know how to act in mathematics, that there is no enigma in that. We regarded, therefore, Kripke's paradox as some kind of cathartic applied to some wild version of this `no mystery' viewpoints\footnote{We have tried to argue that Kripke's challenge acquires some relevance in the context of ifthenism; our use of the skeptic paradox could be considered as enabling some kind of \textit{reductio ad absurdum}.}. We are many times guided to believe that there is an untouchable, mythical and original intuition intended by some formal system, we learn the development of the calculus as something that could not be otherwise. This identification between the necessity of the formal rules in a calculus and the necessity of the calculus itself is what we have labeled as ifthenism\footnote{Note that the Putnam-Finkelstein answer is directed towards posing a challenge to the way we follow formal rules (that is, towards Kripke's challenge against rule-following in general, as our ability to follow rules). What this article is intended to cover is the other half of the problem, what Wittgenstein called the `entanglement' in the rules of a system. Note that, in order to this view to succeed, one must be able to read the quoted paragraphs from PI as stating something different from the powerful and appealing reading that Kripke makes of them.}. What pluralism vindicates is a way of understanding actual mathematical practice, in the sense that every calculus is valid. This naturally leads to some tensions; since we have multiple possibilities of development, one could be inclined to appeal repeatedly to the \textit{original intuition dogma}, and what this shows is that this philosophical picture is untenable: it is one of those we should modify. What we proposed is, in some sense, to make clear that, at each point, we are taking a definite path, that we make some decision. So, in the end, one could argue that we have just changed some picture by another. Well, this is true in some sense, but the \textit{attitude} towards the way mathematics and philosophy meet has changed too. We have seen that there is no reason in believing that these two activities must be separated: in fact, one could say that philosophy should always be carried within mathematics, that philosophy makes explicit some part of mathematical practice on which the mathematician does not shed light. 
\\

To put it in other words, philosophy does not change mathematics, it only changes our approach when concerned with (philosophical) features present in mathematical activity. It would be ridiculous to argue against modal logic on the ground of not bearing -as justification- a direct connection with \textit{reality}. In the same way, it would be absurd to consider some philosophical attitude as false because some theorem could be read as stating something against it (think, as Wittgenstein often did, about some moral law). With this article we hope to have, at least, made the reader wonder about the usual way mathematics and philosophy are linked in the literature: there is no hidden enigma in the way we act, but this does not imply that we do not act in a \textit{certain} way one can point at.

\end{document}